\documentclass[leqno,b5paper]{amsart}

\usepackage{amsmath,amsfonts,amsthm,latexsym}
\usepackage{mathrsfs, textcomp}
\usepackage{tikz}
\usepackage{tikz-cd}
\usetikzlibrary{automata, arrows.meta, positioning, decorations.markings}
\newtheorem{theorem}{Theorem}[section]
\newtheorem{lemma}[theorem]{Lemma}
\newtheorem{proposition}[theorem]{Proposition}
\newtheorem{corollary}[theorem]{Corollary}

\theoremstyle{definition}
\newtheorem{definition}[theorem]{Definition}
\newtheorem{example}[theorem]{Example}
\newtheorem{obs}[theorem]{Observation}
\newtheorem{remark}[theorem]{Remark}
\newtheorem{comment}[theorem]{Comment}

\numberwithin{equation}{section}

\usepackage{amssymb}
\usepackage{amscd}
\usepackage{eqlist, color}
\usepackage{mathrsfs, textcomp}
\usepackage{color, eqlist}
\usepackage{graphicx}
\usepackage{bm}
\usepackage{amsmath,amsfonts,amsthm,latexsym}
\usepackage[all]{xy}

\DeclareMathOperator{\Coz}{Coz}

\DeclareMathOperator{\Rs}{Rs}

\DeclareMathOperator{\Rmt}{Rem}
\DeclareMathOperator{\Nd}{Nd}

\theoremstyle{definition}\newtheorem{thm}{Theorem}[section]
\theoremstyle{definition}
\theoremstyle{definition}
\theoremstyle{definition}
\theoremstyle{definition}
\theoremstyle{definition}
\theoremstyle{remark}
\theoremstyle{definition}

\theoremstyle{definition}
\newtheorem{observation}[thm]{Observation}

\begin{document}
	\title[Some variants of remote sublocales]
	{Some variants of remote sublocales}

	\author{Mbekezeli Nxumalo}
	\address{Department of Mathematical Sciences, University of South Africa, P.O. Box 392, 0003 Pretoria, SOUTH AFRICA.}
	\address{Department of Mathematics, Rhodes University, P.O. Box 94, Grahamstown 6140, South Africa.}
	\email{sibahlezwide@gmail.com}
	\subjclass[2010]{06D22}
	\keywords {remote sublocale, remote point, Booleanization, remote from a dense sublocale, $^{*}$remote from a dense sublocale, $f$-remote preserving, $\gamma$-map}
	\thanks{This paper is part of a Ph.D thesis written under the supervision of Prof. Themba Dube. The author  acknowledges funding from the National Research Foundation of South Africa under Grant 134159}
	\dedicatory{}
	
	
	\let\thefootnote\relax\footnote{}
	
	\begin{abstract}
		We introduce and study some variants of remote sublocales, namely sublocales that are remote from dense sublocales and those that are $^{*}$remote from dense sublocales. We show that the coframe of sublocales coincides with the collection of all sublocales remote from the Booleanization. Furthermore, the supplement of the Booleanization of any locale is the largest sublocale $^{*}$remote from the Booleanization. We give conditions on localic maps such that their induced localic image and pre-image functions preserve sublocales that are remote (resp. $^{*}$remote) from dense sublocales. We introduce new types of localic maps called $f$-remote preserving maps. Given a commuting square \begin{equation}
			\begin{tikzcd}
				{S} \arrow{rr}{g} \arrow{dd}[swap]{\alpha} & & {T} \arrow{dd}{\omega}\\
				& &  \\
				{L} \arrow{rr}[swap]{f} && {M} 
			\end{tikzcd}
		\end{equation} of localic maps, we show that $g$ is $f$-remote preserving precisely when $f[\mathfrak{B}L]$ is remote from $S$.  \textbf{}
	\end{abstract}
	
	\maketitle
	
	
	\section*{Introduction}\label{sect0}
	
	In 1962, Fine and Gillman \cite{FG} defined a \textit{remote point} as a point $p\in \beta \mathbb{R}$ such that $p\notin\overline{N}^{\beta \mathbb{R}}$ for any discrete $N\subseteq\mathbb{R}$. Woods \cite{W}, in 1971, showed that in a metric space $X$, a point $p\in \beta X\smallsetminus X$ is remote if and only if $p\notin\overline{N}^{\beta X}$ for every nowhere dense $N\subseteq X$. The concept of remote point in terms of nowhere dense subsets was further studied by several authors. For instance, Van Douwen \cite{{vD}}, in 1981, gave several characterizations of these points, and Dube \cite{D0}, in 2009, investigated remote points in the category of frames. In 1982, Van Mill \cite{vM} defined a \textit{remote collection} as a collection of closed subsets of a Tychonoff space in which some member of the collection misses every nowhere dense subset of the space. His work motivated the study of remote sublocales in the category of locales that I undertook in 2023, \cite{N}. A sublocale is \textit{remote} in case it misses every nowhere dense sublocale.
	
	In this article, we introduce some types of remote sublocales, namely sublocales that are remote from dense sublocales and those that are $^{*}$remote from dense sublocales. The motivation for these sublocales stems from the fact that remote points in terms of nowhere dense sets initially belonged to the set $\beta X\smallsetminus X$  and were defined with respect to nowhere dense subsets of $X$. 
	
	This article is organized as follows. The first section consists of the neccessary background. Sublocales that are remote from dense sublocales and those that are $^{*}$remote from dense sublocales are introduced in the second section. We characterize them and compare their classes with the class of remote sublocales. The third section focuses on a relationship between these variants of remoteness with the Booleanization. We show that the supplement of the Booleanization of any locale is the largest sublocale $^{*}$remote from the Booleanization. The fourth section considers preservation and reflection of remote (resp. $^{*}$remote) sublocales by localic maps. The work done in this section extends to the fifth section where we discuss $f$-remote preserving and $f$-$^{*}$remote preserving localic maps. We give a condition such that $f$-remote preserving maps are precisely those that preserve remote sublocales.

	
	\section{Preliminaries}\label{sect1}
	\subsection{Locales and Sublocales}\label{subsect11}
	
	For basic locale and frame theory, see \cite{PP1} and \cite{J}. 
	
	We shall use the terms frames and locales interchangeably.

	The \textit{top} element and the \textit{bottom} element of a locale $L$ will be denoted by $1_{L}$ and $0_{L}$, respectively, with subscripts dropped when $L$ is clear from the context. An element $p$ of $L$ is a \textit{point} provided that $p<1$ and $a\wedge b\leq p$ implies that either $a\leq p$ or $b\leq p$, for each $a,b\in L$. The \textit{pseudocomplement} of an element $a\in L$ is denoted by $a^{\ast}$. If $a^{\ast}=0$, then $a$ is \textit{dense} and if $a\vee a^{\ast}=1$, then $a$ is \textit{complemented}. A \textit{Boolean} locale is one in which every element is complemented. The \textit{completely below} relation is denoted by $\prec\!\!\prec$ and a locale is \textit{completely regular} if each of its elements is the join of elements completely below it.

	A \textit{sublocale} of a locale $L$ is a subset $S\subseteq L$ such that (i) $S$ is closed under all meets, and (ii) $x\rightarrow s\in S$ for all $x\in L$ and $s\in S$, where $\rightarrow$ is the Heyting operation on $L$ satisfying: $$a\leq b\rightarrow c\quad\text{if and only if}\quad a\wedge b\leq c$$ for every $a,b,c\in L$. We denote by $\mathcal{S}(L)$ the coframe of sublocales of a locale $L$. The \textit{smallest sublocale} of $L$ is the sublocale $\mathsf{O}=\{1\}$. A sublocale $S$ of $L$ is \textit{void} if $S=\mathsf{O}$, otherwise it is \textit{non-void}.  
	We shall use the prefix $S$- for localic properties defined on a sublocale $S$ of $L$. A sublocale $S\subseteq L$ is \textit{complemented} if it has a complement in $\mathcal{S}(L)$. We denote by $L\smallsetminus S$ the supplement of a locale $S\subseteq L$. 
	The sublocales  $${\mathfrak{c}}(a)=\{x\in L:a\leq x\}\quad\text{and}\quad \mathfrak{o}(a)=\{a\rightarrow x:x\in L\}$$ are, respectively, the \textit{closed} and the \textit{open} sublocales induced by $a\in L$, and are complements of each other. 
	The \textit{closure} of $S\in\mathcal{S}(L)$ is denoted by $\overline{S}$, and $S$ is \textit{dense} if $\overline{S}=L$ (equivalently, if $0_{L}\in S$). The \textit{Booleanization} of $L$ is the sublocale $\mathfrak{B}(L)=\{x\rightarrow 0:x\in L\}$, and is the smallest dense sublocale of $L$. A noteworthy result about dense sublocales is that pseudocomplementation on a dense sublocale is precisely that in the locale. This is so because, if $A$ is a dense sublocale of $L$ and $x\in A$, then writing $x^{*A}$ for the pseudocomplement in $x$ in $A$, we have the equalities
	$$x^{*A}=x\rightarrow_{A} 0_{A} =x\rightarrow 0_{L}=x^{*}.$$
	This means $\mathfrak{B}S=\mathfrak{B}L$ for a dense sublocale $S$ of $L$. \cite{P} A sublocale $S$ of $L$ is \textit{nowhere dense} if $S\cap\mathfrak{B}L=\mathsf{O}$. By Dube \cite{D1}, for each $a\in L$, $\mathfrak{c}(a)$ is nowhere dense if and only if $a$ is dense. More generally, a sublocale $N$ of a locale $L$ is nowhere dense if and only if $\bigwedge\!N$ is dense in $L$, \cite{N}.
	
	A \textit{localic map} is an infima preserving map between locales such that the corresponding left adjoints preserve finite meets. Associated with every localic map is its left adjoint called the \textit{frame homomorphism} which preserves binary meets and arbitrary joins. We shall write $f:L\rightarrow M$ for a localic map and $f^{*}$ for its left adjoint. A frame homomorphism $f^{*}:M\rightarrow L$ is called \textit{dense} in case it maps only the bottom element to the bottom element; a \textit{quotient map} provided that it is surjective; an \textit{extension} if it is a dense quotient map; and \textit{closed} if $f(x\vee f^{*}(y))=f(x)\vee y$ for all $x\in L$ and $y\in M$.  According to Dube \cite{D2}, a \textit{nowhere dense} map is a frame homomorphism $f^{*}:M\rightarrow L$ in which for each non-zero $x\in M$ there is a non-zero $y\in M$ with $y\leq x$ such that $f^{*}(y)=0_{L}$. Every sublocale $S\subseteq L$ is associated with the quotient map $\nu_{S}:L\rightarrow S$ defined by $$\nu_{S}(a)={\bigwedge}\{s\in S: a\leq s\}.$$
	
	For each $a\in L$ and $S\in\mathcal{S}(L)$, $\mathfrak{o}_{S}(\nu_{S}(a))=S\cap \mathfrak{o}(a)$ and ${\mathfrak{c}_{S}}(\nu_{S}(a))=S\cap \mathfrak{c}(a).$
	
	Each localic map $f:L\rightarrow M$ induces the maps $$f[-]:\mathcal{S}(L)\rightarrow \mathcal{S}(M)\text{ and }f_{-1}[-]:\mathcal{S}(M)\rightarrow\mathcal{S}(L)$$ which are called the \textit{localic image function} induced by $f$ and the \textit{localic preimage function} induced by $f$, respectively. The maps $f[-]$ and $f_{-1}[-]$ form an adjunction.
	
		\subsection{On $\beta L$, $\lambda L$ and $\upsilon L$} Refer to \cite{BM,MV,DN1,BG1} for more details on $\beta L$, $\lambda L$ and $\upsilon L$. We represent by $\beta L$ the \textit{Stone-\v{C}ech compactification} of a completely regular locale $L$. The frame homomorphism $\beta_{L}:\beta L\rightarrow L$, defined by $I\mapsto \bigvee\!I$, is an extension with its right adjoint denoted by $r_{L}$. For any localic map $f:L\rightarrow M$ between completely regular locales $L$ and $M$, the left adjoint $\beta(f)^{*}$ of $\beta(f)$ is given by $$\beta(f)^{*}: J\mapsto \{x\in L: x\leq f^{*}(y)\text{ for some }y\in J\}.$$

	We use $\Coz L$ to denote the \textit{cozero part} of $L$. The \textit{regular Lindelöf reflection} of $L$ is the locale of $\sigma$-ideals of $\Coz L$ and is denoted by ${\lambda}L$. The join map ${\lambda}_{L}:{\lambda}L\rightarrow L$ is an extension. We define the extension $k_{L}:\beta L\rightarrow\lambda L$ by $I\mapsto\langle I\rangle_{\sigma}$, where $\langle\cdot\rangle$ signifies $\sigma$-ideal generation in $\Coz L$. 
	
	For a completely regular locale $L$, the \textit{realcompact reflection} of $L$ is the locale $\upsilon L$ defined to be Fix($\ell$), where  \[\ell:{\lambda}L\rightarrow{\lambda}L, I\mapsto\Big[{\bigvee} I\Big]\wedge{\bigwedge}\{P\in Pt({\lambda} L): I\leq P\}.\]The join map $v_{L}:vL\rightarrow L$ is an extension and $\ell_{L}:\lambda L\rightarrow \upsilon L$ is an extension effected by $\ell$.
	
	When $\beta L$ is regarded as the locale of regular ideals of $\Coz L$, we get the following commuting diagram in the category of completely regular locales whose morphisms are localic maps between them.
	\begin{equation}\label{CRegFrmdiagram}
		\begin{gathered}
			\xymatrixcolsep{4pc}\xymatrixrowsep{2.5pc}
			\xymatrix{
				{L} \ar[rd]^{(\upsilon_{L})_{*}}  \ar[rdd]_{(\lambda_{L})_{*}} \ar[ddd]_{r_{L}} \ar[rrr]^{f}&& &{M} \ar[dl]_{(\upsilon_{M})_{*}} \ar[ddd]^{r_{M}}\ar[ddl]^{(\lambda_{L})_{*}}\\
				&\upsilon L\ar[r]^{\upsilon(f)} \ar[d]^{(\ell_{L})_{*}}&\upsilon M \ar[d]_{(\ell_{M})_{*}} & \\
				&\lambda L\ar[r]_{\lambda(f)} \ar[dl]^{(\kappa_{L})_{*}}&\lambda M \ar[dr]_{(\kappa_{M})_{*}} & \\
				\beta L \ar[rrr]_{\beta f} &&& \beta M}.
		\end{gathered}
	\end{equation}By a \textit{$\gamma$-lift} we mean the localic morphism $\gamma(f):\gamma L\rightarrow \gamma M$, where $\gamma\in \{\beta,\lambda,\upsilon\}$. 


	\section{Remoteness from a dense sublocale}
	This section focuses on some variants of remoteness which are defined with respect to dense sublocales. 
	
	We shall write $\overline{N}$ for the closure of a sublocale $N$ of a locale $L$, taken in $L$ and write $\overline{N}^{S}$ for the closure of a sublocale $N$ taken in a sublocale $S$ of a locale $L$. In the case of open and closed sublocales of a locale $L$, we shall write $\mathfrak{c}(a)$ (or $\mathfrak{c}_{L}(a)$) and $\mathfrak{o}(a)$ (or $\mathfrak{o}_{L}(a)$), respectively. To refer to a closed and open sublocale of a sublocale $S$ of $L$, we shall, respectively, write  $\mathfrak{c}_{S}(a)$ and $\mathfrak{o}_{S}(a)$.

	We remind the reader that prefix $S$- will be used to indicate a localic property defined in a sublocale $S$ of a locale $L$. For instance, if $S$ is a sublocale of a locale $L$ and $N\in\mathcal{S}(S)$, then $S$-nowhere dense $N$ means that $N$ is nowhere dense in a sublocale $S$ of $L$.
\begin{definition}\label{remote}
	Let $S\subseteq L$ be a dense sublocale of $L$. Then
	\begin{enumerate}
		\item  $T\in \mathcal{S}(L)$ is \textit{remote from $S$} if $T\cap \overline{N}=\mathsf{O}$ for every $S$-nowhere dense $N\in\mathcal{S}(S)$.
		\item  A sublocale $T\subseteq L\smallsetminus S$ is \textit{{*}remote from $S$} if $T\cap \overline{N}=\mathsf{O}$ for every $S$-nowhere dense $N\in\mathcal{S}(S)$.
	\end{enumerate}
\end{definition} 


	We note some examples.
\begin{example}\label{BLandL}
	(1) $\mathsf{O}$ is remote from every dense sublocale of a locale $L$.

	(2) Recall from \cite{D0} that a point $I$ of $\beta L$ is \textit{remote} provided that $I\vee r(h_{*}(0))=\top$ for every nowhere dense quotient map $h:L\rightarrow M$. It is easy to see that a point $I\in\beta L$ is remote if and only if $\mathfrak{c}(I)$ is remote from $L$.

	(3) In \textbf{Top}, we say that $A\subseteq X\smallsetminus Y$, where $Y$ is a dense subspace of $X$, is \textit{$^{*}$remote from $Y$} in case $A\cap \overline{N}^{X}=\emptyset$ for all nowhere dense subsets $N$ of $Y$. For a Tychonoff space $X$, a point $p\in\beta X\smallsetminus X$ is remote if and only if $\{p\}$ is $^{*}$remote from $X$. 
	
	(4) For each $A,B\in\mathcal{S}(L)$, if $A\subseteq B$ and $B$ is remote from $S$, then $A$ is remote from $S$.
\end{example}
	
	Denote by $\mathcal{S}_{\text{rem}}(L\ltimes S)$ and 
	$^{*}\mathcal{S}_{\text{rem}}(L\ltimes S)$ the collections of sublocales that are remote and $^{*}$remote from a dense sublocale $S$, respectively.


	
	

	We characterize members of $\mathcal{S}_{\text{rem}}(L\ltimes S)$ and ${^{*}\mathcal{S}_{\text{rem}}(L\ltimes S)}$. Recall that for each $S\in\mathcal{S}(L)$ and any $a\in S$,
$\overline{\mathfrak{c}_{S}(a)}= \mathfrak c\Big(\bigwedge(\mathfrak c_S(a))\Big)= \mathfrak c(a).$

\begin{theorem}\label{opendensefrom}
	Let $S\in\mathcal{S}(L)$ be dense and $A\in\mathcal{S}(L)$ $ ( $resp. $A\in\mathcal{S}(L\smallsetminus S)$$ ) $. The following statements are equivalent.
	\begin{enumerate}
		\item $A\in\mathcal{S}_{\text{rem}}(L\ltimes S)$ $ ( $resp. $A\in{^{*}\mathcal{S}_{\text{rem}}(L\ltimes S)}$$)$.
		\item For all $S$-dense $x\in S$, $A\cap\mathfrak{c}(x)=\mathsf{O}$.
		\item For all $S$-dense $x\in S$, $A\subseteq\mathfrak{o}(x)$.
		\item For every $S$-dense $x\in S$, $\nu_{A}(x)=1$.
	\end{enumerate}
\end{theorem}
\begin{proof}
	$(1)\Rightarrow(2)$: Let $x\in S$ be $S$-dense. Then $\mathfrak{c}_{S}(x)$ is $S$-nowhere dense. By (1), $$\mathsf{O}=A\cap \overline{\mathfrak{c}_{S}(x)}=A\cap\mathfrak{c}(x).$$
	
	$(2)\Rightarrow(3)$: Follows since for all $x\in S$, $A\cap \mathfrak{c}(x)=\mathsf{O}$ if and only if $A\subseteq \mathfrak{o}(x)$.
	
	$(3)\Rightarrow(4)$: Follows since for any $x\in S$, $\nu_{A}(x)=1$ if and only if $A\subseteq \mathfrak{o}(x)$. 
	
	$(4)\Rightarrow(1)$: Let $N\in\mathcal{S}(S)$ be $S$-nowhere dense. Then $\bigwedge\!N$ is $S$-dense. By (4), $\nu_{A}\left(\bigwedge\!N\right)=1$ so that $$\mathsf{O}=\mathfrak{c}_{A}\left(\nu_{A}\left(\bigwedge\!N\right)\right)=A\cap\mathfrak{c}\left(\bigwedge N\right)=A\cap\overline{N}$$as required.
\end{proof}
\begin{comment}
The equivalence of (4) in Theorem \ref{opendensefrom} tells us that, $A$ is remote (resp. *remote) from $S$ if and only if it is contained in every open dense sublocale of $L$ induced by an element of $S$. This is reminiscent of the characterization of the remote sublocales of $L$ as precisely those that are contained in every dense sublocale of $L$.
\end{comment}

The preceding theorem leads us to the following example of a sublocale of $L$ which is remote from a dense sublocale $S$. 

\begin{example}\label{ND(S)remotefrom}
	Set $$\Nd(S)=\bigvee\{N\in\mathcal{S}(L):N\text{ is nowhere dense}\}.$$ If $S$ is a dense sublocale of $L$, then the sublocale $L\smallsetminus \overline{\Nd(S)}$ is remote from $S$. To see this, choose an $S$-dense $x\in S$, then $\mathfrak{c}_{S}(x)$ is $S$-nowhere dense and contained in $\Nd(S)$, so that $\overline{\mathfrak{c}_{S}(x)}=\mathfrak{c}(x)\subseteq\overline{\Nd(S)}$. Therefore $\mathfrak{c}(x)\cap (L\smallsetminus \overline{\Nd(S)}) =\mathsf{O}$. By Theorem \ref{opendensefrom}(2), $L\smallsetminus\overline{\Nd(S)}$ is remote from $S$.

	Observe that $L\smallsetminus \overline{\Nd(S)}$ may be different from $\mathsf{O}$. Consider a dense sublocale $S\in\mathcal{S}(L)$ where $\Nd(S)$ is $S$-nowhere dense and $L\neq\mathsf{O}$ (for instance, a locale whose Booleanization is complemented, see \cite[Corollary 4.16]{DS}). Since $\Nd(S)$ is the largest $S$-nowhere dense sublocale and its closure in $S$ is $S$-nowhere dense, $\Nd(S)=\overline{\Nd(S)}^{S}$ making it $S$-closed nowhere dense. Because $S$ is dense in $L$, $\Nd(S)$ is nowhere dense in $L$ so that $\overline{\Nd(S)}$ is nowhere dense in $L$. Therefore $L\neq \overline{\Nd(S)}$ which means that $L\smallsetminus\overline{\Nd(S)}\neq \mathsf{O}$.
\end{example}

Since the set $\mathcal{S}_{\text{rem}}(L\ltimes S)$ does not restrict where its members come from, we have the following immediate connection between $\mathcal{S}_{\text{rem}}(L\ltimes S)$ and 
$^{*}\mathcal{S}_{\text{rem}}(L\ltimes S)$.
\begin{proposition}\label{remotesets}
	For every dense sublocale $S$ of a locale $L$, $^{*}\mathcal{S}_{\text{rem}}(L\ltimes S)\subseteq\mathcal{S}_{\text{rem}}(L\ltimes S)$.
\end{proposition}

\begin{obs}
	For a non-void Boolean locale $L$, we have that $^{*}\mathcal{S}_{\text{rem}}(L\ltimes L)\subset\mathcal{S}_{\text{rem}}(L\ltimes L)$. This is because $L\in\mathcal{S}_{\text{rem}}(L\ltimes L)$ but $L$ is not contained in $L\smallsetminus L=\mathsf{O}$.
\end{obs}
	
\begin{remark}
In an attempt to obtain an equality in Proposition \ref{remotesets}, we start by recalling from \cite{P} that a sublocale is \textit{rare} if its supplement is the whole locale. Restricting our sublocales to dense and rare sublocales yields the following:

	\begin{quote}
	\emph{If $S$ is a dense and rare sublocale of a locale $L$, then $^{*}\mathcal{S}_{\text{rem}}(L\ltimes S)=\mathcal{S}_{\text{rem}}(L\ltimes S)$.}
\end{quote}
Sublocales which are simultaneously dense and rare do exist. For instance, recall that Plewe in \cite{P} defines a locale to be \textit{dense in itself} if every Boolean sublocale has a dense supplement. He then shows that a locale is dense in itself if and only if its Booleanization is rare. The locale $\mathfrak{O}(\mathbb{R})$, where $\mathbb{R}$ is the set of real numbers, is an example of a dense in itself locale. This follows from the fact that the space $\mathbb{R}$ is dense in itself (because it has no isolated points) and since, according to \cite{P}, every sober space is dense in itself if and only if its locale of opens is dense in itself, the real line being sober and dense in itself makes $\mathfrak{O}(\mathbb{R})$ dense in itself.

\end{remark}
We have the following relationship between the collection $\mathcal{S}_{\text{rem}}(L)$ of all remote sublocales of a locale $L$ and $\mathcal{S}_{\text{rem}}(L\ltimes S)$.

\begin{proposition}\label{SRemandSRemLS}
	Let $L$ be a locale. Then $\mathcal{S}_{\text{rem}}(L)\subseteq \mathcal{S}_{\text{rem}}(L\ltimes S)$ for every dense sublocale $S$ of $L$. 
\end{proposition} 
\begin{proof}
	Follows since a remote sublocale of $L$ misses the closure of every nowhere dense sublocale of $L$, and hence the closure of every nowhere dense sublocale of $S$ since $S$ is a dense sublocale of $L$. 
\end{proof}


	The following result shows a relationship between the collections $\mathcal{S}_{\text{rem}}(S)$ and $\mathcal{S}_{\text{rem}}(L\ltimes S)$.
\begin{proposition}\label{remS}
	Let $S$ be a dense sublocale of a locale $L$. Then \[\mathcal{S}(S)\cap \mathcal{S}_{\text{rem}}(L\ltimes S)=\mathcal{S}_{\text{rem}}(S).\]
\end{proposition}
\begin{proof}
	$\mathcal{S}(S)\cap \mathcal{S}_{\text{rem}}(L\ltimes S)\subseteq\mathcal{S}_{\text{rem}}(S)$: Let $A\in \mathcal{S}(S)\cap \mathcal{S}_{\text{rem}}(L\ltimes S)$ and choose an $S$-nowhere dense $N\in\mathcal{S}(S)$. Then $A\cap \overline{ N}=\mathsf{O}$ which implies that $A\cap N=\mathsf{O}$. Thus $A\in \mathcal{S}_{\text{rem}}(S)$.

	$\mathcal{S}_{\text{rem}}(S)\subseteq\mathcal{S}(S)\cap \mathcal{S}_{\text{rem}}(L\ltimes S)$: Let $A\in \mathcal{S}_{\text{rem}}(S)$ and choose an $S$-nowhere dense $N$. Then $\overline{N}^{S}$ is $S$-nowhere dense so that $\mathsf{O}=A\cap \overline{N}^{S}=A\cap \overline{N}\cap S=A\cap\overline{N}$. Thus $A\in \mathcal{S}(S)\cap \mathcal{S}_{\text{rem}}(L\ltimes S)$.
\end{proof}
	
	We close this section with a discussion of elements inducing closed sublocales that are remote and $^{*}$remote from dense sublocales.
	
	Set
\[\Rmt(L\ltimes S )=\{a\in L:\mathfrak c_L(a)\in\mathcal{S}_{\text{rem}}(L\ltimes S)\}\]and
\[^{*}\Rmt(L\ltimes S )=\{a\in L:\mathfrak c_L(a)\in {^{*}\mathcal{S}_{\text{rem}}(L\ltimes S)}\}.\]

We give the following proposition, the proof of which follows easily from the definitions.
\begin{proposition}\label{sublocale}
	Let $S\in\mathcal{S}(L)$ be dense. 
	\begin{enumerate}
		\item $a\in \Rmt(L\ltimes S )$ iff $a\vee x=1$ for all dense $x\in S$.
		\item For $\mathfrak{c}(a)\subseteq L\smallsetminus S$, $a\in {^{*}\Rmt(L\ltimes S)}$ iff $a\vee x=1$ for all dense $x\in S$.
	\end{enumerate}
\end{proposition}


\section{Connection between the Booleanization and Remoteness from dense sublocales }

We begin this section with a result showing that the Booleanization of a locale is remote from every dense sublocale.

\begin{proposition}\label{BLisremote}
	Let $L$ be a locale. Then $\mathfrak{B}L$ is remote from every dense sublocale of $L$.
\end{proposition}
\begin{proof}
	This follows since $\mathfrak{B}L$ is a remote sublocale of $L$ so that, by Proposition \ref{SRemandSRemLS}, it is remote from every dense sublocale of $L$.
\end{proof}
The above result tells us that $\mathfrak{B}L\in \mathcal{S}_{\text{rem}}(L\ltimes\mathfrak{B}L)$. 
	In what follows we show that $\mathcal{S}_{\text{rem}}(L\ltimes\mathfrak{B}L)$ is in fact the whole coframe $\mathcal{S}(L)$.  Observe that if $K$ is a dense sublocale of $L$, then  $\mathfrak{B}L=\mathfrak{B}K$, as a consequence of which a sublocale of $K$ is nowhere dense in $K$ if and only if it is nowhere dense in $L$.

\begin{proposition}\label{rempropBL}
	Let $L$ be a locale. Then  $\mathcal{S}_{\text{rem}}(L\ltimes\mathfrak{B}L)=\mathcal{S}(L)$. 
\end{proposition}
\begin{proof}
	If $A\in \mathcal{S}(L)$ and $N\in\mathcal{S}(\mathfrak{B}L)$ is nowhere dense in $\mathfrak{B}L$, then $N=\mathsf{O}$ which implies that $A\cap \overline{N}=\mathsf{O}$. Thus $A\in \mathcal{S}_{\text{rem}}(L\ltimes\mathfrak{B}L)$. Hence $\mathcal{S}(L)\subseteq\mathcal{S}_{\text{rem}}(L\ltimes\mathfrak{B}L)$, making $\mathcal{S}_{\text{rem}}(L\ltimes\mathfrak{B}L)=\mathcal{S}(L)$ since the other containment always holds. 
\end{proof}

\begin{obs}\label{obsremotefrom}
	Using Proposition \ref{rempropBL} and the fact that $L\in\mathcal{S}_{\text{rem}}(L)$ if and only if $L$ is Boolean (from \cite[Proposition 3.5.]{N}), it is easy to see that for a non-Boolean locale $L$, $L\in\mathcal{S}_{\text{rem}}(L\ltimes\mathfrak{B}L)$ but $L\notin\mathcal{S}_{\text{rem}}(L)$. This is a particular case where we do not have equality in Proposition \ref{SRemandSRemLS}. However, $\mathcal{S}_{\text{rem}}(L\ltimes L)=\mathcal{S}_{\text{rem}}(L)$.
\end{obs}





It turns out that the Booleanization is the only sublocale remote from itself, as we show below. 

\begin{theorem}
	Let $S$ be a dense sublocale of $L$. The following statements are equivalent.
	\begin{enumerate}
		\item $S\in\mathcal{S}_{\text{rem}}(L\ltimes S)$. 
		\item $S=\mathfrak{B}L$.
		\item $L$ is remote from $S$.
	\end{enumerate}
\end{theorem}
\begin{proof}
	$(1)\Rightarrow(2)$: Let $N$ be $S$-nowhere dense. Then since $S\in\mathcal{S}_{\text{rem}}(L\ltimes S)$, $S\cap \overline{N}=\mathsf{O}$ which implies that $\mathsf{O}=S\cap N=N$. This makes $S$ Boolean. So, $S=\mathfrak{B}L$ because the only dense Boolean sublocale of $L$ is $\mathfrak{B}L$.
	
	
	$(2)\Rightarrow (3)$: Since $\mathsf{O}$ is the only nowhere dense sublocale of $\mathfrak{B}L$, we have that $L\cap \overline{N}=\mathsf{O}$ for every nowhere dense sublocale $N$ of $\mathfrak{B}L=S$. 

	$(3)\Rightarrow(1)$: Let $N$ be a nowhere dense sublocale of $S$. Since $L$ is remote from $S$, $L\cap\overline{N}=\mathsf{O}$, making $N=\mathsf{O}$. Hence $S\cap\overline{N}=\mathsf{O}$ so that $S\in\mathcal{S}_{\text{rem}}(L\ltimes S)$. 
\end{proof}
	
	We showed in \cite{N} that the Booleanization is the largest remote sublocale of a locale. We work towards establishing a relationship between the Booleanization and the largest sublocale that is remote from dense sublocales. Using Theorem \ref{opendensefrom}(2), one can easily see that the join of remote sublocales is remote. This tells us that every locale has the largest sublocale which is remote from a specific dense sublocale. For $S\in\mathcal{S}(L)$ dense in $L$, set $$\Rs(L\ltimes S)=\bigvee\{A\in \mathcal{S}(L):A\text{ is remote from }S\}.$$ 
	\begin{remark}\label{Lislarge}
		Proposition \ref{rempropBL} tells us that $L=\Rs(L\ltimes \mathfrak{B}L)$.
	\end{remark}
	
We give the following lemma to prepare for Theorem \ref{RsBL}.
	
	\begin{lemma}\label{SRemLemma}
		Let $L$ be a locale. If $A\in\mathcal{S}_{\text{rem}}(L\ltimes S)$, then $A\cap S\in\mathcal{S}_{\text{rem}}(L)$. 
	\end{lemma} 
	\begin{proof}
		Assume that $A\in\mathcal{S}_{\text{rem}}(L\ltimes S)$ and let $N$ be nowhere dense in $L$. Since $N\cap S\subseteq N$, $N\cap S$ is nowhere dense in $L$ which in turns makes it $S$-nowhere dense. By hypothesis, $A\cap \overline{S\cap N}=\mathsf{O}$. This makes $(A\cap S)\cap N=\mathsf{O}$. Thus $A\cap S\in\mathcal{S}_{\text{rem}}(L)$.
	\end{proof}
	
	\begin{theorem}\label{RsBL}
		For a dense sublocale $S$ of a locale $L$, $\Rs(L\ltimes S)\cap S=\mathfrak{B}L$.
	\end{theorem}

	\begin{proof}
		 $\mathfrak{B}L\subseteq \Rs(L\ltimes S)\cap S$: Follows since $\mathfrak{B}L$ is remote from $S$ and $\mathfrak{B}L\subseteq S$ by density of $S$. 
		 
		 
		 For the other containment, using Lemma \ref{SRemLemma} we get that $\Rs(L\ltimes S)\cap S\in \mathcal{S}_{\text{rem}}(L)$. $\mathfrak{B}L$ being the largest remote sublocale of $L$ gives $\Rs(L\ltimes S)\cap S\subseteq\mathfrak{B}L$. Thus $\Rs(L\ltimes S)\cap S=\mathfrak{B}L$.
	\end{proof}
\begin{observation}\label{BLnotremote}
	The statement ``$\Rs(L\ltimes S)=\mathfrak{B}L$ for every dense sublocale $S$ of $L$" is not true for a non-Boolean locale $L$. Otherwise, $\mathfrak{B}L=\Rs(L\ltimes\mathfrak{B}L)=L$ where the latter equality follows from Remark \ref{Lislarge}, which is not possible.
\end{observation}

We noticed in Example \ref{ND(S)remotefrom} that $L\smallsetminus\overline{\Nd(S)}$ is remote from a dense sublocale $S$ of $L$. A question about a relationship between the sublocales $L\smallsetminus\overline{\Nd(S)}$ and $\Rs(L\ltimes S)$ arises. We address this in the following result. 

\begin{proposition}
	Let $S$ be a dense sublocale of a locale $L$. The following statements are equivalent:
	\begin{enumerate}
		\item $\Rs(L\ltimes S)=L\smallsetminus\overline{\Nd(S)}$. 
		\item $\Nd(S)$ is $S$-nowhere dense.
	\end{enumerate} 
\end{proposition}
\begin{proof}
	$(1)\Rightarrow(2)$: Assume that $\Rs(L\ltimes S)=L\smallsetminus\overline{\Nd(S)}$. Since $\mathfrak{B}L\subseteq \Rs(L\ltimes S)$, we have that $\mathfrak{B}L\subseteq L\smallsetminus\overline{\Nd(S)}$. Therefore $\mathfrak{B}L\cap \overline{\Nd(S)}=\mathsf{O}$ which implies that $\mathfrak{B}L\cap\Nd(S)=\mathsf{O}$, making $\bigwedge\Nd(S)$ dense in $L$. But $S$ is dense, so $\nu_{S}\left(\bigwedge\Nd(S)\right)=\bigwedge\Nd(S)$ is $S$-dense. It follows that $\Nd(S)$ is $S$-nowhere dense.

	$(2)\Rightarrow(1)$: We show that $\Rs(L\ltimes S)\subseteq L\smallsetminus\overline{\Nd(S)}$. $\Nd(S)$ being $S$-nowhere dense implies that $\Rs(L\ltimes S)\cap \overline{\Nd(S)}=\mathsf{O}$. This gives $\Rs(L\ltimes S)\subseteq L\smallsetminus\overline{\Nd(S)}$ as required.
\end{proof}

We noticed in Observation \ref{BLnotremote} that for a non-Boolean locale $L$, $\mathfrak{B}L\neq\Rs(L\ltimes \mathfrak{B}L)$. For $^{*}$remoteness, we show in Theorem \ref{RsDense} that the supplement of $\mathfrak{B}L$, for any locale $L$, is the largest sublocale $^{*}$remote from $\mathfrak{B}L$.

Just like in the case of $\Rs(L\ltimes S)$, set $$^{*}\!\Rs(L\ltimes S)=\bigvee\{A\in \mathcal{S}(L):A\text{ is }^{*}\text{remote from }S\}$$ for a dense sublocale $S$ of $L$.

We give the following result which we shall use below. 
\begin{lemma}\label{rempropBLstar}
	$^{*}\mathcal{S}_{\text{rem}}(L\ltimes\mathfrak{B}L)=\{T\in\mathcal{S}(L):T\subseteq L\smallsetminus\mathfrak{B}L\}$ for every locale $L$.
\end{lemma}
\begin{obs}\label{obsremotefromstar}
	(1) From Lemma \ref{rempropBLstar}, observe that when $L$ is not dense in itself, we get another case where $^{*}\mathcal{S}_{\text{rem}}(L\ltimes S)\neq\mathcal{S}_{\text{rem}}(L\ltimes S)$. This is because we have that $L\neq L\smallsetminus\mathfrak{B}L$ so that by Proposition \ref{rempropBL}, $L\in \mathcal{S}_{\text{rem}}(L\ltimes\mathfrak{B}L)$ but $L\notin{^{*}\mathcal{S}_{\text{rem}}(L\ltimes\mathfrak{B}L)}$.
	
	(2) A locale is dense in itself if and only if it is $^{*}$remote from its Booleanization: Observe that $L$ is dense in itself if and only if $\mathfrak{B}L$ is rare if and only if $L\subseteq L\smallsetminus\mathfrak{B}L$ if and only if $L\in\{T\in\mathcal{S}(L):T\subseteq L\smallsetminus\mathfrak{B}L\}={^{*}\mathcal{S}_{\text{rem}}(L\ltimes\mathfrak{B}L)}$, where the last equality holds by Lemma \ref{rempropBLstar}.
\end{obs}

		\begin{theorem}\label{RsDense}
		Let $L$ be a locale. Then ${^{*}\!\Rs(L\ltimes \mathfrak{B}L)}=L\smallsetminus\mathfrak{B}L$.
	\end{theorem}
	\begin{proof}Follows since, by Lemma \ref{rempropBLstar}, $L\smallsetminus\mathfrak{B}L\in {^{*}\mathcal{S}_{\text{rem}}(L\ltimes\mathfrak{B}L)}$, making $L\smallsetminus\mathfrak{B}L\subseteq{^{*}\!\Rs(L\ltimes \mathfrak{B}L)}$.  Also all ${^{*}}$remote sublocales (including ${^{*}\Rs(L\ltimes \mathfrak{B}L)}$) belong to the set $\{T\in\mathcal{S}(L):T\subseteq L\smallsetminus\mathfrak{B}L\}$. 
	\end{proof} 
	
	
		

	\section{Preservation and reflection of sublocales that are remote (resp. *remote) from dense sublocales}
	
	In this section we discuss localic maps that send back and forth the sublocales introduced in Definition \ref{remote}. 
	
	Consider a commuting diagram  
	\begin{equation}\label{famousdia}
		\begin{tikzcd}
			{S} \arrow{rr}{g} \arrow{dd}[swap]{\alpha} & & {T} \arrow{dd}{\omega}\\
			& &  \\
			{L} \arrow{rr}[swap]{f} && {M} 
		\end{tikzcd}
	\end{equation}
	where $S,T,L$ and $M$ are locales, $f$ and $g$ are localic maps and the downward morphisms are dense injective localic maps. Our discussion will make use of the information provided in diagram \ref{famousdia}. We commented in the preliminaries that a localic map $k:P\rightarrow Q$ is dense if and only if $k[P]$ is a dense sublocale of $Q$. So, $\alpha[S]$ and $\omega[T]$ are dense sublocales of $L$ and $M$, respectively. Since for a quotient map $v:W\rightarrow Y$, $v_{*}:Y\rightarrow v_{*}[Y]$ is an isomorphism, we will sometimes write $S$ and $T$ for the sublocales $\alpha[S]$ and $\omega[T]$, respectively. 
	
	Diagram \ref{CRegFrmdiagram} is a particular case of the situation depicted in diagram \ref{famousdia}.

	For $^{*}\!$remoteness, we note that, $A\cap \alpha[S]=\mathsf{O}$ implies $A\subseteq L\smallsetminus \alpha[S]$ but $A\subseteq L\smallsetminus \alpha[S]$ does not always imply that $A$ misses $\alpha[S]$ unless $\alpha[S]$ is complemented. We will sometimes treat these cases differently. 
	

	
	We start by recording a description of localic maps that preserve remoteness and $^{*}$remoteness from dense sublocales. For the following result, we recall from \cite{FPP} that a localic map $f:L\rightarrow M$ \textit{takes $A$-remainder to $B$-remainder} if $f[L\smallsetminus A]\subseteq M\smallsetminus B$ where $A\in\mathcal{S}(L)$, $B\in\mathcal{S}(M)$. We shall write $f:L\rightarrow M$ takes $S$-remainder to $T$-remainder  to mean that $f$ takes $\alpha[S]$-remainder to $\omega[T]$-remainder. Recall from \cite{M3} that a frame homomorphism $f^{*}:M\rightarrow L$ is said to be \textit{weakly closed} in case $a\vee f^{*}(b)=1$ implies $f(a)\vee b=1$ for every $a\in L$ and $b\in M$. 
	
	\begin{proposition}\label{beta}
		Assume that $g^{*}$ in diagram \ref{famousdia} is skeletal  and $f^{*}\circ \omega=\alpha\circ g^{*}$. Then
		\begin{enumerate}
			\item $f[{\mathcal{S}_{\text{rem}}(L\ltimes S)}]\subseteq {\mathcal{S}_{\text{rem}}( M\ltimes T)}$.
			\item If $f^{*}$ is weakly closed, then $f[{\Rmt( L\ltimes S)}]\subseteq{\Rmt(M\ltimes T)}$.
		\end{enumerate}
	\end{proposition}
	\begin{proof}
		(1) Let $A\in {\mathcal{S}_{\text{rem}}(L\ltimes S)}$ and choose an $\omega[T]$-dense $y\in \omega[T]$. Then $y=\omega(t)$ for some $t\in T$. Because $g^{*}$ is skeletal, we have that $g^{*}(t)$ is $S$-dense so that $\alpha(g^{*}(t))$ is $\alpha[S]$-dense since $\alpha:S\rightarrow\alpha[S]$ is an isomorphism. Therefore $\mathsf{O}=A\cap \mathfrak{c}(\alpha(g^{*}(t))=A\cap \mathfrak{c}(f^{*}(\omega(t))$ where the latter equality follows since $f^{*}\circ \omega=\alpha\circ g^{*}$. Therefore $$A\subseteq \mathfrak{o}(f^{*}(\omega(t)))=f_{-1}[\mathfrak{o}(\omega(t))]=f_{-1}[\mathfrak{o}(y)].$$ We get that $f[A]\subseteq f[f_{-1}[\mathfrak{o}(y)]]\subseteq\mathfrak{o}(y)$. This tells us that $f[A]$ is contained in every open sublocale induced by an $\omega[T]$-dense element, so by Theorem \ref{opendensefrom}(3), $f[A]\in{\mathcal{S}_{\text{rem}}( M\ltimes T)} $.

		(2) Let $x\in \Rmt(L\ltimes S)$ and choose an $\omega[T]$-dense $t\in \omega[T]$. Therefore $t=\omega(y)$ for some $y\in T$. Since $g^{*}$ is skeletal, $g^{*}(y)$ is $S$-dense, making $\alpha(g^{*}(y))$ $\alpha[S]$-dense. It follows that $\alpha(g^{*}(y))\vee x=1$. Because $\alpha\circ g^{*}=f^{*}\circ \omega$, $f^{*}(\omega(y))\vee x=1$. The weakly closedness of $f^{*}$ implies that $1_{M}=\omega(y)\vee f(x)=t\vee f(x)$. Thus $f(x)\in\Rmt(M\ltimes T)$.
	\end{proof}
	
	\begin{obs}
		In terms of $\gamma$-lifts, the condition $f^{*}\circ \omega=\alpha\circ g^{*}$ on $f$ resembles that of a \textit{$\gamma$-map} which was defined in \cite{DN1} as a frame homomorphism $t:M\rightarrow L$ that satisfies $\gamma(t)\circ (\gamma_{M})_{\ast}=(\gamma_{L})_{\ast}\circ t.$
	\end{obs}
\begin{proposition}\label{betastar}
	If the map $g^{*}$ in diagram \ref{famousdia} is skeletal, $f^{*}\circ \omega=\alpha\circ g^{*}$ and $f$ takes $S$-remainder to $T$-remainder, then:
		\begin{enumerate}
			\item $f[{^{*}\!\mathcal{S}_{\text{rem}}(L\ltimes S)}]\subseteq {^{*}\!\mathcal{S}_{\text{rem}}( M\ltimes T)}$.
			\item If $f^{*}$ is weakly closed, then $f[{^{*}\!\Rmt( L\ltimes S)}]\subseteq{^{*}\!\Rmt(M\ltimes T)}$.
		\end{enumerate}
\end{proposition}
\begin{proof}
			With the assumption that $f$ takes $S$-remainder to $T$-remainder, it is clear that $A\subseteq L\smallsetminus \alpha[S]$ implies $f[A]\subseteq f[L\smallsetminus \alpha[S]]\subseteq M\smallsetminus \omega[T]$ for all $A\in\mathcal{S}(L)$. This together with Proposition \ref{beta}(1)\&(2) show that both (1) and (2) hold.
\end{proof}

	\begin{obs}
	For Proposition \ref{betastar}, in the case where $\alpha[S]$ is complemented, we replace $f$ takes $S$-remainder to $T$-remainder with the condition that $f$ is injective and $f[\alpha[S]]=\omega[T]$. From this we get that $A\subseteq L\smallsetminus \alpha[S]$ implies $A\cap \alpha[S]=\mathsf{O}$. Therefore $\mathsf{O}=f[\mathsf{O}]=f[A\cap \alpha[S]]=f[A]\cap f[\alpha[S]]$ so that $f[A]\subseteq M\smallsetminus f[\alpha[S]]=M\smallsetminus \omega[T]$.
\end{obs}
	

	
	We return to descriptions of localic maps that reflect and preserve the variants of remoteness introduced in Definition \ref{remote}.
	
	\begin{proposition}\label{beta1}
		Assume that the morphism $g$ in diagram \ref{famousdia} is skeletal. Then
		\begin{enumerate}
			\item $f[A]\in {\mathcal{S}_{\text{rem}}(M\ltimes T)}$ implies $A\in {\mathcal{S}_{\text{rem}}(L\ltimes S)}$ for every $A\in\mathcal{S}(L)$.
			\item $f(x)\in \Rmt(M\ltimes T)$ implies $x\in \Rmt(L\ltimes S)$ for all $x\in L$.
		\end{enumerate}
	\end{proposition}
	\begin{proof}
		(1) Assume that $f[A]\in {\mathcal{S}_{\text{rem}}(M\ltimes T)}$ and let $a\in \alpha[S]$ be $\alpha[S]$-dense. Then $a=\alpha(x)$ for some $x\in S$ where such $x$ is $S$-dense. Since $g$ is skeletal, $g(x)$ is $T$-dense so that $\omega(g(x))$ is $\omega[T]$-dense. It follows that $f[A]\subseteq \mathfrak{o}(\omega(g(x)))$ which implies $f[A]\subseteq\mathfrak{o}(f(\alpha(x)))$ because $k\circ g=f\circ \alpha$. Therefore $$A\subseteq f_{-1}[f[A]]\subseteq f_{-1}[\mathfrak{o}(f(\alpha(x)))]= \mathfrak{o}(f^{*}(f(\alpha(x))))\subseteq\mathfrak{o}(\alpha(x))=\mathfrak{o}(a).$$ Thus $A\in {\mathcal{S}_{\text{rem}}(L\ltimes S)}$.

		(2) Follows similar sketch of the proof of (1).
	\end{proof}
		\begin{proposition}\label{beta1star}
		If $g$ in diagram \ref{famousdia} is skeletal,  $\omega[T]$ is a complemented sublocale of $M$ and $f_{-1}[\omega[T]]=\alpha[S]$, then:
			\begin{enumerate}
				\item $f[A]\in {^{*}\!\mathcal{S}_{\text{rem}}(M\ltimes T)}$ implies $A\in {^{*}\!\mathcal{S}_{\text{rem}}(L\ltimes S)}$ for every $A\in\mathcal{S}(L)$.
				\item $f(x)\in {^{*}\!\Rmt(M\ltimes T)}$ implies $x\in {^{*}\!\Rmt(L\ltimes S)}$ for all $x\in L$.
			\end{enumerate}
	\end{proposition}
	\begin{proof}
We only show that $f[A]\subseteq M\smallsetminus \omega[T]$ implies $A\subseteq L\smallsetminus \alpha[S]$. Observe that for complemented $\omega[T]$ in $M$ with $f_{-1}[\omega[T]]=\alpha[S]$,
	\begin{align*}
		f[A]\subseteq M\smallsetminus \omega[T]
		&\quad\Longleftrightarrow\quad
		f[A]\cap \omega[T]=\mathsf{O}\\
		&\quad\Longrightarrow\quad
		A\cap f_{-1}[\omega[T]]=\mathsf{O}\\
		&\quad\Longleftrightarrow\quad
		A\cap \alpha[S]=\mathsf{O}\\
		&\quad\Longrightarrow\quad
		A\subseteq L\smallsetminus \alpha[S].
	\end{align*} for all $A\in\mathcal{S}(L)$. 
\end{proof}
	\begin{proposition}\label{for}
		Suppose that the localic map $g$ in diagram \ref{famousdia} is skeletal. Then:
		\begin{enumerate}
			\item $f_{-1}[{\mathcal{S}_{\text{rem}}(M\ltimes T)}]\subseteq{\mathcal{S}_{\text{rem}}(L\ltimes S)}$.
			\item $f^{*}[{\Rmt(M\ltimes T)}]\subseteq{\Rmt(L\ltimes S)}$.
		\end{enumerate}
	\end{proposition}
	\begin{proof}
		(1) Let $A\in {\mathcal{S}_{\text{rem}}(M\ltimes T)}$ and choose an $S$-dense $a\in S$. Then $a=\alpha(x)$ for some $x\in S$ which is $S$-dense. $g(x)$ is $T$-dense because $g$ is skeletal. It follows that $A\subseteq \mathfrak{o}(\omega(g(x)))$ since $\omega(g(x))$ is $\omega[T]$-dense and $A\in {\mathcal{S}_{\text{rem}}(M\ltimes T)}$. Therefore \[f_{-1}[A]\subseteq f_{-1}[\mathfrak{o}(\omega(g(x)))]=\mathfrak{o}(f^{*}(\omega(g(x))))=\mathfrak{o}(f^{*}(f(\alpha(x))))\subseteq\mathfrak{o}(\alpha(x))=\mathfrak{o}(x)\] making $f_{-1}[A]\in\mathcal{S}_{\text{rem}}(L\ltimes S)$.
		
		(2) Proof follows similar sketch of the proof of (1).
	\end{proof}

	\begin{proposition}\label{forstar}
	Suppose that the localic map $g$ in diagram \ref{famousdia} is skeletal, $f_{-1}[\omega[T]]=\alpha[S]$ and $\omega[T]$ is complemented in $M$, then:
		\begin{enumerate}
			\item $f_{-1}[^{*}{\mathcal{S}_{\text{rem}}(M\ltimes T)}]\subseteq{^{*}\!\mathcal{S}_{\text{rem}}(L\ltimes S)}$.
			\item $f^{*}[^{*}{\Rmt(M\ltimes T)}]\subseteq{^{*}\!\Rmt(L\ltimes S)}$.
		\end{enumerate} 
\end{proposition}
		\begin{proof}Assume that $f_{-1}[\omega[T]]=\alpha[S]$ and $\omega[T]$ is complemented in $M$. We only show that $A\subseteq M\smallsetminus \omega[T]$ implies $f_{-1}[A]\subseteq L\smallsetminus \alpha[S]$ which is needed for both (1) and (2). We have that $A\subseteq M\smallsetminus \omega[T]$ gives $A\cap \omega[T]=\mathsf{O}$. Therefore $\mathsf{O}=f_{-1}[A]\cap f_{-1}[\omega[T]]=f_{-1}[A]\cap \alpha[S]$, which implies that $f_{-1}[A]\subseteq L\smallsetminus \alpha[S]$.
	\end{proof}
	\begin{proposition}\label{for1}
		If $g^{*}$ in diagram \ref{famousdia} is skeletal, $\alpha\circ g^{*}=f^{*}\circ \omega$ and $f[-]$ is surjective, then:
			\begin{enumerate}
				\item $f_{-1}[A]\in {\mathcal{S}_{\text{rem}}(L\ltimes S)}$ implies $A\in {\mathcal{S}_{\text{rem}}(M\ltimes T)}$ for all $A\in\mathcal{S}(M)$.
				\item If $f$ takes $S$-remainder to $T$-remainder, then $f_{-1}[A]\in {^{*}\!\mathcal{S}_{\text{rem}}(L\ltimes S)}$ implies $A\in {^{*}\!\mathcal{S}_{\text{rem}}(M\ltimes T)}$ for every $A\in\mathcal{S}(M)$.
			\end{enumerate}
	\end{proposition}
	\begin{proof}
		(1) Let $A\in\mathcal{S}(M)$ be such that $f_{-1}[A]\in\mathcal{S}_{\text{rem}}(L\ltimes S)$ and choose an $\omega[T]$-dense $b\in \omega[T]$. Then $b=\omega(x)$ for some $x\in T$ with $x$ a $T$-dense element. The skeletalness of $g^{*}$ implies that $g^{*}(x)$ is $S$-dense so that $\alpha(g^{*}(x))$ is $\alpha[S]$-dense. Therefore $$\mathsf{O}=f_{-1}[A]\cap \mathfrak{c}(\alpha(g^{*}(x)))=f_{-1}[A]\cap \mathfrak{c}(f^{*}(\omega(x)))=f_{-1}[A]\cap f_{-1}[\mathfrak{c}(\omega(x))]=f_{-1}[A\cap\mathfrak{c}(\omega(x))].$$Since $f[-]$ is surjective,  $\mathsf{O}=f[f_{-1}[A\cap\mathfrak{c}(\omega(x))]]=A\cap \mathfrak{c}(\omega(x))=A\cap \mathfrak{c}(b)$. Thus $A\in\mathcal{S}_{\text{rem}}(M\ltimes T)$.

		(2) We only show that $f_{-1}[A]\subseteq L\smallsetminus \alpha[S]$ implies $A\subseteq L\smallsetminus \alpha[S]$. Observe that $$
		f_{-1}[A]\subseteq L\smallsetminus \alpha[S]\;\Longrightarrow\;
		f[f_{-1}[A]]\subseteq f[L\smallsetminus \alpha[S]]\;\Longrightarrow\;
		A\subseteq M\smallsetminus \omega[T]$$ which proves the result.
	\end{proof}
	\begin{proposition}\label{for1star}
	Assume that $g^{*}$ in diagram \ref{famousdia} is skeletal and one of the following statements holds:
	\begin{enumerate}
		\item[(a)] $f^{*}$ is weakly closed  $g$ is surjective.
		\item[(b)] $\alpha\circ g^{*}=f^{*}\circ \omega$ and $f$ is surjective.
	\end{enumerate} Then
		\begin{enumerate}
			\item For each $x\in M$, $f^{*}(x)\in {\Rmt(L\ltimes S)}$ implies $x\in {\Rmt(M\ltimes T)}$.
			\item If $f[L\smallsetminus \alpha[S]]\subseteq M\smallsetminus \omega[T]$, then $f^{*}(x)\in {^{*}\!\Rmt(L\ltimes S)}$ implies $x\in {^{*}\!\Rmt(M\ltimes T)}$ for all $x\in M$.
		\end{enumerate}
\end{proposition}
	\begin{proof}
	(1) Suppose that $f^{*}$ is weakly closed, $g$ is surjective and let $T\in T$ be $T$-dense. Then $t=\omega(b)$ for some $b\in T$ which is $T$-dense. Then $g^{*}(b)$ is $S$-dense since $g^{*}$ is skeletal. Because $f^{*}(x)\in\Rmt(L\ltimes S)$, we get that \begin{equation}\label{fclosed}
		f^{*}(x)\vee \alpha(g^{*}(b))=1.
	\end{equation} The weakly closedness of $f^{*}$ gives $x\vee f(\alpha(g^{*}(b)))=1$. Therefore $$1=x\vee \omega(g(g^{*}(b)))=x\vee \omega(b)=x\vee t$$ where the second equality holds since $g$ is surjective. Thus $x\in \Rmt(M\ltimes T)$.

	Assume that $\alpha\circ g^{*}=f^{*}\circ \omega$ and $f$ is surjective. Then from  equation \ref{fclosed}, we get that $f^{*}(x)\vee f^{*}(t)=1$ which implies that $$f(f^{*}(x\vee t))=f(f^{*}(x)\vee f^{*}(t))=f(1)=1$$ so that by surjectivity of $f$, $x\vee t=1$ making $x\in \Rmt(M\ltimes T)$.

	(2) Can be deduced from Proposition \ref{for1}(2) and (1) above.
\end{proof}


\section{$f$-remote preserving and $f$-$^{*}\!$remote preserving maps}

In this section, we pay a closer attention to localic maps with the properties given in Proposition \ref{beta}(1) and Proposition \ref{betastar}(1). We still make use of diagram \ref{famousdia}.

We give the following definition. 

\begin{definition}
	We call a map $g$ in diagram \ref{famousdia} \textit{$f$-remote preserving} if $f[\mathcal{S}_{\text{rem}}(L\ltimes S)]\subseteq\mathcal{S}_{\text{rem}}(M\ltimes T)$ and \textit{$f$-$^{*}\!$remote preserving} if $f[^{*}\!\mathcal{S}_{\text{rem}}(L\ltimes S)]\subseteq{^{*}\!\mathcal{S}_{\text{rem}}(M\ltimes T)}$.  
\end{definition}

Since $\mathfrak{B}L\in \mathcal{S}_{\text{rem}}(L\ltimes S)$ for every dense sublocale $S$ of $L$, in the next result, we characterize $f$-remote preserving maps in terms of the Booleanization of a locale. We also include, in the same result, a characterization in terms of the largest sublocale remote from a given dense sublocale. We recall that if $w:P\rightarrow Q$ is a dense injective localic map, then for all $x\in P$, $x$ is $P$-dense if and only if $w(x)$ is $L$-dense.

\begin{theorem}\label{gammaremotepreserving}
	Suppose that $f^{*}\circ \omega=\alpha\circ g^{*}$. The following statements are equivalent.
	\begin{enumerate}
		\item $g$ is $f$-remote preserving.
		\item $f[\mathfrak{B}L]\in\mathcal{S}_{\text{rem}}(M\ltimes T)$.
		\item $f[\Rs(L\ltimes S)]\in \mathcal{S}_{\text{rem}}(M\ltimes T)$.
		\item $f[\Rs( L\ltimes S)]\subseteq\Rs(M\ltimes T)$.
	\end{enumerate}
\end{theorem}
\begin{proof}
	$(1)\Longrightarrow (2)$: Since $\mathfrak{B}L$ is remote from every dense sublocale of $L$ and $\alpha[S]$ is a dense sublocale of $L$, $\mathfrak{B}L$ is remote from $\alpha[S]$. By (1), $f[\mathfrak{B}L]\in \mathcal{S}_{\text{rem}}(M\ltimes T)$.

	$(2)\Longrightarrow(3)$: Let $a\in \omega[T]$ be $\omega[T]$-dense. Then $a=\omega(x)$ for some $x\in T$. By hypothesis, $f[\mathfrak{B} L]\subseteq\mathfrak{o}(\omega(x))$. Therefore $\mathfrak{B} L\subseteq f_{-1}[\mathfrak{o}(\omega(x))]=\mathfrak{o}[f^{*}(\omega(x))]$. Since $f^{*}\circ \omega=\alpha\circ g^{*}$, $\mathfrak{B} L\subseteq\mathfrak{o}[\alpha(g^{*}(x))]$, making $\alpha(g^{*}(x))$ $L$-dense so that $g^{*}(x)$ is $S$-dense and hence $\alpha(g^{*}(x))$ is $\alpha[S]$-dense. But $\Rs(L\ltimes S)$ is remote from $\alpha[S]$, so $\Rs(L\ltimes S)\subseteq\mathfrak{o}[\alpha(g^{*}(x))]$. Therefore \[f[\Rs(L\ltimes S)]\subseteq f[\mathfrak{o}(\alpha(g^{*}(x)))]=f[\mathfrak{o}(f^{*}(\omega(x)))]=f[f_{-1}[\mathfrak{o}(\omega(x))]]\subseteq\mathfrak{o}(\omega(x))=\mathfrak{o}(a).\]Thus $f[\Rs(L\ltimes S)]\in \mathcal{S}_{\text{rem}}(M\ltimes T)$.

	$(3)\Longrightarrow(4)$: Follows since $\Rs(M\ltimes T)$ is the largest sublocale remote from $\omega[T]$.

	$(4)\Longrightarrow (1)$: Let $A\in\mathcal{S}_{\text{rem}}(L\ltimes S)$. Then $A\subseteq\Rs(L\ltimes S)$ so that $f[A]\subseteq f[\Rs(L\ltimes S)]$. But $f[\Rs(L\ltimes S)]\subseteq\Rs(M\ltimes T)$, so $f[A]\subseteq f[\Rs(M\ltimes T)]$. Since sublocales contained in members of $\mathcal{S}_{\text{rem}}(M\ltimes T)$ are remote from $\omega[T]$, $f[A]$ is remote from $\omega[T]$. Thus $g$ is $f$-remote preserving.
\end{proof}
We give the following characterization of $f$-$^{*}\!$remote preserving maps.
\begin{proposition}\label{stargammaremotepreserving}
	Assume that $f^{*}\circ \omega=\alpha\circ g^{*}$. The following statements are equivalent.
	\begin{enumerate}
		\item $g$ is $f$-$^{*}\!$remote preserving.
		\item $f[^{*}\!\Rs(L\ltimes S)]$ is $^{*}\!$remote from $\omega[T]$.
		\item $f[^{*}\!\Rs( L\ltimes S)]\subseteq{^{*}\!\Rs(M\ltimes T)}$.
	\end{enumerate}
\end{proposition}
\begin{proof}
	$(1)\Longrightarrow(2)$: Follows since ${^{*}\!\Rs(L\ltimes S)}\in{^{*}\!\mathcal{S}_{\text{rem}}(L\ltimes S)}$.
	
	$(2)\Longrightarrow(3)$: Trivial.
	
	$(3)\Longrightarrow (1)$: Proof is similar to that of Theorem \ref{gammaremotepreserving} $(4)\Longrightarrow(1)$.
\end{proof}

In Proposition \ref{remotepreservation} below, we explore a relationship between $f$-remote preservation and preservation of remote sublocales. We give the following lemma which will be useful in proving the result.

\begin{lemma}\label{gammapreservationlemma}
	The following statements hold.
	\begin{enumerate}
		\item $A\in\mathcal{S}_{\text{rem}}(S)$ iff $\alpha[A]\in\mathcal{S}_{\text{rem}}(L\ltimes S)$.
		\item $A\in \mathcal{S}_{\text{rem}}(L\ltimes S)$ implies $\alpha_{-1}[A]\in\mathcal{S}_{\text{rem}}(S)$.
	\end{enumerate}
\end{lemma}
\begin{proof}
	(1) Recall from Proposition \ref{remS} that  $\mathcal{S}_{\text{rem}}(K)=\mathcal{S}_{\text{rem}}(L\ltimes K)\cap \mathcal{S}(K)$ for each dense $K\in\mathcal{S}(L)$. Therefore $\alpha[A]\in\mathcal{S}_{\text{rem}}(\alpha[S])$ if and only if $\alpha[A]\in\mathcal{S}_{\text{rem}}(L\ltimes S)\cap\mathcal{S}(\alpha[S])$. Since $\alpha:S\rightarrow\alpha[S]$ is an isomorphism, it is easy to see that $A\in\mathcal{S}_{\text{rem}}(S)$ if and only if $\alpha[A]\in\mathcal{S}_{\text{rem}}(\alpha[S])$ if and only if $\alpha[A]\in\mathcal{S}_{\text{rem}}(L\ltimes S)\cap\mathcal{S}(\alpha[S])$.

	(2) Let $x\in S$ be $S$-dense. Then $\alpha(x)$ is $\alpha[S]$-dense. It follows that $A\cap \mathfrak{c}_{L}(\alpha(x))=\mathsf{O}$. Therefore $$\mathsf{O}=\alpha_{-1}[A]\cap \alpha_{-1}[\mathfrak{c}_{L}(\alpha(x))]=\alpha_{-1}[A]\cap \mathfrak{c}_{S}((\alpha)^{*}(\alpha(x)))=\alpha_{-1}[A] \cap \mathfrak{c}_{S}(x)$$proving the result.
\end{proof}
Recall from \cite[Theorem 4.1.]{N} that a localic map $f:L\rightarrow M$ preserves remote sublocales if and only if $f[\mathfrak{B}L]$ is a remote sublocale of $M$.
\begin{proposition}\label{remotepreservation}
	Assume that $f^{*}\circ \omega=\alpha\circ g^{*}$. Then $g$ is $f$-remote preserving iff $g[-]$ preserves remote sublocales.
\end{proposition}
\begin{proof}
	$(\Longrightarrow):$ Since $\alpha[S]$ is dense in $L$, $\mathfrak{B}\alpha[S]=\mathfrak{B}L$ making $\mathfrak{B}\alpha[S]=\alpha[\mathfrak{B}S]$ remote from $\alpha[S]$. Because $g$ is $f$-remote preserving, we have that $f[\alpha[\mathfrak{B}S]]\in\mathcal{S}_{\text{rem}}(M\ltimes T)$ which implies that $\omega[g[\mathfrak{B}S]]\in\mathcal{S}_{\text{rem}}(M\ltimes T)$ since $\omega^{*}\circ f=g\circ \alpha^{*}$. It follows from Lemma \ref{gammapreservationlemma}(2) that $(\omega)_{-1}[\omega[g[\mathfrak{B}S]]]\in \mathcal{S}_{\text{rem}}(T)$. But $g[\mathfrak{B}S]\subseteq(\omega)_{-1}[\omega[g[\mathfrak{B}S]]]$, so $g[\mathfrak{B}S]\in\mathcal{S}_{\text{rem}}(T)$. By \cite[Theorem 4.1.]{N}, $g$ preserves remote sublocales.

	$(\Longleftarrow):$ We show that $f[\mathfrak{B}L]$ is remore from $\omega[T]$. Since $\mathfrak{B}L\in\mathcal{S}_{\text{rem}}(L\ltimes S)$, it follows from Lemma \ref{gammapreservationlemma}(2) that $(\alpha)_{-1}[\mathfrak{B}L]\in\mathcal{S}_{\text{rem}}(S)$. By hypothesis, $g[(\alpha)_{-1}[\mathfrak{B}L]]\in\mathcal{S}_{\text{rem}}(T)$. By Lemma \ref{gammapreservationlemma}(1), $\omega[g[(\alpha)_{-1}[\mathfrak{B}L]]]\in\mathcal{S}_{\text{rem}}(M\ltimes T)$ which implies that $f[\alpha[(\alpha)_{-1}[\mathfrak{B}L]]]\in\mathcal{S}_{\text{rem}}(M\ltimes T)$ since $f\circ \alpha=\omega\circ g$. But $\mathfrak{B}L=\mathfrak{B}\alpha[S]\subseteq \alpha[S]$ and using the fact that $\alpha: S\rightarrow \alpha[S]$ is an isomorphism, $$f[\alpha[(\alpha)_{-1}[\mathfrak{B}L]]]=f[\alpha[(\alpha)_{-1}[\mathfrak{B}\alpha[S]]]]=f[\mathfrak{B}\alpha[S]]=f[\mathfrak{B}L]\in\mathcal{S}_{\text{rem}}(M\ltimes T)$$ as required.
\end{proof}




Consider a commuting diagram 
\begin{equation}\label{STRULM}
	\begin{tikzcd}
		{S} \arrow{rd}{i} \arrow{dd}[swap]{\alpha} \arrow{rrr}{g}&& & {T} \arrow{dl}[swap]{k} \arrow{dd}{\omega}&\\
		&R\arrow{r}{\varphi} \arrow{dl}{\theta}&U \arrow{dr}[swap]{\sigma} & && & \\
		{L} \arrow{rrr}[swap]{f} && & {M} &
	\end{tikzcd}
\end{equation}
where $S,T,R,U,L$ and $M$ are locales, the downward arrows are dense injective localic maps and the horizontal arrows are localic maps. We find a relationship between $f$-remote preservation and $\varphi$-remote preservation. En route to that, we give the following lemma.

\begin{lemma}\label{bvl} From diagram \ref{STRULM}, $\theta[\mathcal{S}_{\text{rem}}(R\ltimes S)]\subseteq\mathcal{S}_{\text{rem}}(L\ltimes S)$.
\end{lemma}
\begin{proof}
	We have that $\alpha$ is dense since it is the composite of two dense localic maps $i$ and $\theta$. Let $A\in\mathcal{S}_{\text{rem}}(R\ltimes S)$ and choose an $\alpha[S]$-dense $y\in \alpha[S]$. Then $y=\alpha(x)$ for some $x\in S$. Since $\alpha:S\rightarrow \alpha[S]$ is an isomorphism, $x$ is $S$-dense so that $i(x)$ is $i[S]$-dense. Therefore $A\cap \mathfrak{c}_{R}(i(x))=\mathsf{O}$. Observe that $\theta[A]\cap \mathfrak{c}_{L}(\alpha(x))=\mathsf{O}$. To see this, let $a\in \theta[A]\cap \mathfrak{c}_{L}(\alpha(x))$. Then $a=\theta(b)$ for some $b\in A$ and $\alpha(x)\leq a$. We have that $$i(x)=\theta^{*}(\theta(i(x)))=\theta^{*}(\alpha(x))\leq \theta^{*}(\theta(b))=b$$ since $\theta$ is injective and $\alpha=\theta\circ i$.  Therefore $b\in A\cap \mathfrak{c}_{R}(i(x))$ which implies $b=1$ so that $a=\theta(b)=1$. Thus $\theta[A]\cap \mathfrak{c}_{L}(\alpha(x))=\mathsf{O}$. Hence $\theta[A]\in \mathcal{S}_{\text{rem}}(L\ltimes S)$.
\end{proof}
Since $\theta[R]\subseteq L$, we have that \[\{B\in \mathcal{S}(\theta[R]):B\cap \alpha[S]=\mathsf{O}\}\subseteq\{C\in \mathcal{S}(L):C\cap \alpha[S]=\mathsf{O}\}\] so that \[\theta[R]\smallsetminus \alpha[S]=\bigvee\{B\in \mathcal{S}(\theta[R]):B\cap \alpha[S]=\mathsf{O}\}\subseteq\bigvee\{C\in \mathcal{S}(L):C\cap \alpha[S]=\mathsf{O}\}=L\smallsetminus \alpha[S].\] As a result of this and Lemma \ref{bvl}, we have the following result.
\begin{lemma}\label{starbvl} From diagram \ref{STRULM}, $\theta[{^{*}\!\mathcal{S}_{\text{rem}}(R\ltimes S)}]\subseteq{^{*}\!\mathcal{S}_{\text{rem}}(L\ltimes S)}$.
\end{lemma}
\begin{obs}\label{bvlobs}
	In light of the preceding two lemmas and the relationship between the $\beta,\upsilon$ and $\lambda$ extensions depicted in diagram \ref{CRegFrmdiagram}, we have \[\mathcal{S}_{\text{rem}}(\upsilon L\ltimes L)\subseteq\mathcal{S}_{\text{rem}}(\lambda L\ltimes L)\subseteq\mathcal{S}_{\text{rem}}(\beta L\ltimes L)\] and \[{^{*}\!\mathcal{S}_{\text{rem}}(\upsilon L\ltimes L)}\subseteq{^{*}\!\mathcal{S}_{\text{rem}}(\lambda L\ltimes L)}\subseteq{^{*}\!\mathcal{S}_{\text{rem}}(\beta L\ltimes L)}.\]
\end{obs}
\begin{proposition}\label{gfremote}
	If $g$ in diagram \ref{STRULM} is $f$-remote preserving, then it is $\varphi$-remote preserving. 
\end{proposition}
\begin{proof}
	Let $A\in\mathcal{S}_{\text{rem}}(R\ltimes S)$. It follows from Lemma \ref{bvl} that $\theta[A]\in\mathcal{S}_{\text{rem}}(L\ltimes S)$. Since $g$ is $f$-remote preserving, $f[\theta[A]]\in\mathcal{S}_{\text{rem}}(M\ltimes T)$ making $\sigma[\varphi[A]]\in\mathcal{S}_{\text{rem}}(M\ltimes T)$. By Lemma \ref{gammapreservationlemma}(2), $(\sigma)_{-1}[\sigma[\varphi[A]]]\in\mathcal{S}_{\text{rem}}(U)\subseteq \mathcal{S}_{\text{rem}}(U\ltimes T)$. Since $\varphi[A]\subseteq(\sigma)_{-1}[\sigma[\varphi[A]]]$ and a sublocale of any member of $\mathcal{S}_{\text{rem}}(U\ltimes T)$ belongs to $\mathcal{S}_{\text{rem}}(U\ltimes T)$, we have that $\varphi[A]\in \mathcal{S}_{\text{rem}}(U\ltimes T)$. Thus $g$ is $\varphi$-remote preserving.
\end{proof}
\begin{obs}\label{starobsgfremote}
	Recall from \cite{FPP} that given a localic map $f:L\rightarrow M$ and any $K\in\mathcal{S}(L)$, $f[L\smallsetminus K]\subseteq M\smallsetminus f[K]$ whenever $K=f_{-1}[J]$ for some $J\in \mathcal{S}(M)$. Since for the $^{*}\!$remoteness case of Proposition \ref{gfremote} we need $\varphi[R\smallsetminus i[S]]\subseteq U\smallsetminus k[T]$, we assume that  $i[S]=\varphi_{-1}[k[T]]$ and $\varphi[-]$ is surjective in diagram \ref{STRULM}. Then $$\varphi[R\smallsetminus i[S]]\subseteq U\smallsetminus\varphi[i[S]]=U\smallsetminus \varphi[\varphi_{-1}[k[T]]]=U\smallsetminus k[T]$$ so that $A\in {^{*}\!\mathcal{S}_{\text{rem}}}(R\ltimes S)$ implies $\varphi[A]\in {^{*}\!\mathcal{S}_{\text{rem}}}(U\ltimes T)$. This approach also helps in verifying $^{*}\!$remoteness cases of Corollary \ref{blambdavremote} and Proposition \ref{tfg}(2)\&(3) below.
\end{obs}
\begin{obs}\label{obsfremote}
	The converse of Proposition \ref{gfremote} holds if $\alpha[-]$ is surjective (hence an isomorphism). Indeed, assume that $g$ is $\varphi$-remote preserving and let $A\in \mathcal{S}_{\text{rem}}(L\ltimes S)$. By Lemma \ref{gammapreservationlemma}(2), $\alpha_{-1}[A]\in \mathcal{S}_{\text{rem}}(S)$. Therefore $i[\alpha_{-1}[A]]\in \mathcal{S}_{\text{rem}}(R\ltimes S)$
	by Lemma \ref{gammapreservationlemma}(1). Since $g$ is $\varphi$-remote preserving, $\varphi[i[\alpha_{-1}[A]]]\in \mathcal{S}_{\text{rem}}(U\ltimes T)$. By Lemma \ref{bvl}, $\sigma[\varphi[i[\alpha_{-1}[A]]]]\in \mathcal{S}_{\text{rem}}(M\ltimes T)$ so that $f[\theta[i[\alpha_{-1}[A]]]]\in \mathcal{S}_{\text{rem}}(M\ltimes T)$ because $f\circ \theta=\sigma\circ \varphi$. Since $\theta\circ i=\alpha$ and $\alpha[-]$ is surjective, $f[\theta[i[\alpha_{-1}[A]]]]=f[\alpha[\alpha_{-1}[A]]]=f[A]\in \mathcal{S}_{\text{rem}}(M\ltimes T)$. Thus $g$ is $f$-remote preserving.
\end{obs}
Call a localic map $f:L\rightarrow M$ \textit{$\gamma$-remote preserving} if $\gamma(f)\left[\mathcal{S}_{\text{rem}}(\gamma L\ltimes L)\right]\subseteq {\mathcal{S}_{\text{rem}}(\gamma M\ltimes M)}$ and \textit{$\gamma$-$^{*}\!$remote preserving} provided that $\gamma(f)\left[^{*}\!\mathcal{S}_{\text{rem}}(\gamma L\ltimes L)\right]\subseteq {^{*}\!\mathcal{S}_{\text{rem}}(\gamma M\ltimes M)}.$
\begin{corollary}\label{blambdavremote}
	We have that $$\beta\text{-remote preserving}\;\Longrightarrow\;\lambda\text{-remote preserving}\;\Longrightarrow\;\upsilon\text{-remote preserving}.$$
\end{corollary}
\begin{obs}\label{obsbvl}
	To get the reverse directions of Corollary \ref{blambdavremote}, we observe that the morphisms $\upsilon_{L}$, $\beta_{L}$ and $\lambda_{L}$ are isomorphisms whenever $L$ is compact. The case of $\lambda_{L}$ follows since, according to \cite{J0}, $\lambda_{L}$ is injective (hence an isomorphism) whenever $L$ is Lindelöf. Because every compact locale is Lindelöf, we have that $\lambda_{L}$ is an isomorphism whenever $L$ is compact.
\end{obs}
We end this section with the following result. 
\begin{proposition}\label{tfg}
	Consider a commuting diagram
	\begin{equation}
		\begin{tikzcd}
			{L } \arrow{rr}{f}\arrow{ddrr}[swap]{t}& & {M } \arrow{dd}{\varphi} &\\
			& & & {\text{  }} & & & \\
			& & {N} &
		\end{tikzcd}
	\end{equation}where $f,g$ and $t$ are localic maps and $L,M$ and $N$ are locales. \begin{enumerate}
		\item If $\varphi$ and $f$ are $\gamma$-remote preserving, then $t$ is $\gamma$-remote preserving.
		\item If $t$ is $\gamma$-remote preserving, $\varphi$  sends elements to dense elements, then $f$ is $\gamma$-remote preserving.
		\item If $t$ is $\gamma$-remote preserving and $A\subseteq\gamma(f)[\mathfrak{B}_{\gamma L}]$ for all $A\in\mathcal{S}_{\text{rem}}(\gamma M\ltimes M)$, then $\varphi$ is $\gamma$-remote preserving.
	\end{enumerate}
\end{proposition}
\begin{proof}
	(1) For each $A\in\mathcal{S}(\gamma L)$, we have 
	\begin{align*}
		A\in\mathcal{S}_{\text{rem}}(\gamma L\ltimes L)
		&\quad\Longrightarrow\quad
		\gamma(f)[A]\in\mathcal{S}_{\text{rem}}(\gamma M\ltimes M)\quad\text{since}\quad f\quad\text{is}\quad\gamma\text{-remote preserving}\\
		&\quad\Longrightarrow\quad
		\gamma(\varphi)[\gamma(f)[A]]\in\mathcal{S}_{\text{rem}}(\gamma N\ltimes N)\text{ since } \varphi\text{ is }\gamma\text{-remote preserving}\\
		&\quad\Longrightarrow\quad
		\gamma(\varphi\circ f)[A]\in\mathcal{S}_{\text{rem}}(\gamma N\ltimes N)\quad\text{since}\quad \gamma\quad\text{is a functor}\\
		&\quad\Longrightarrow\quad
		\gamma(t)[A]\in\mathcal{S}_{\text{rem}}(\gamma N\ltimes N).
	\end{align*}

	(2) Let $A\in\mathcal{S}_{\text{rem}}(\gamma L\ltimes L)$ and choose dense $x\in M$. Since $t$ is $\gamma$-remote preserving, we have that $\gamma(t)[A]\in\mathcal{S}_{\text{rem}}(\gamma N\ltimes N)$. Since $\varphi$ is skeletal, we have that $\varphi(x)$ is dense in $N$. It follows that $\gamma(t)[A]\subseteq\mathfrak{o}((\gamma_{N})_{*}(\varphi(x)))$. But $(\gamma_{N})_{*}\circ \varphi=\gamma(\varphi)\circ (\gamma_{M})_{*}$, so $\gamma(t)[A]\subseteq\mathfrak{o}(\gamma(\varphi)((\gamma_{M})_{*}(x)))$. Therefore $\gamma(\varphi)_{-1}[\gamma(t)[A]]\subseteq \mathfrak{o}((\gamma_{M})_{*}(x))$ making $\gamma(\varphi)_{-1}[\gamma(t)[A]]$ remote from $M$. Since $\gamma(\varphi)[\gamma(f)[A]]=\gamma(t)[A]$ implies $\gamma(f)[A]\subseteq \gamma(\varphi)_{-1}[\gamma(t)[A]],$ we have that $\gamma(f)[A]$ is remote from $M$.

	(3) Let $A\in\mathcal{S}_{\text{rem}}(\gamma M\ltimes M)$. Then $A\subseteq \gamma(f)[\mathfrak{B}_{\gamma L}]$ which implies that $$\gamma(\varphi)[A]\subseteq\gamma(\varphi)[\gamma(f)[\mathfrak{B}_{\gamma L}]]=\gamma(t)[\mathfrak{B}_{\gamma L}].$$ But $\gamma(t)[\mathfrak{B}_{\gamma L}]$ is remote from $N$, so $\gamma(\varphi)[A]$ is remote from $N$. It follows that $\varphi$ is $\gamma$-remote preserving.
\end{proof}
\begin{obs}
	We commented in Observation \ref{starobsgfremote} that the approach given in that observation can be used to prove the $^{*}\!$remoteness case of Proposition \ref{tfg}(2)\&(3). The $^{*}\!$remoteness case of Proposition \ref{tfg}(1) follows the similar sketch of the proof of Proposition \ref{tfg}(1) where  $\gamma$-remote preserving is replaced by  $\gamma$-${^{*}}$remote preserving.
\end{obs}
	
	

\end{document}